\documentclass{pspum-l}

\newtheorem{theorem}{Theorem}[section]
\newtheorem{lemma}[theorem]{Lemma}
\newtheorem{proposition}[theorem]{Proposition}

\newtheorem{definition}[theorem]{Definition}

\theoremstyle{remark}

\numberwithin{equation}{section}

\usepackage{hyperref,amsmath,mathabx,stackrel}
\usepackage[all]{xy}

\usepackage{tikz-cd}
\usetikzlibrary{arrows.meta}

\newcommand{\opcit}{{\it op.cit.\/}\ }
\newcommand{\ie}{{\it i.e.\/}\ }
\newcommand{\cf}{{\it cf.}}

\def\A{{\mathbb A}}
\def\B{{\mathbb B}}
\def\C{{\mathbb C}}
\def\F{{\mathbb F}}
\def\N{{\mathbb N}}
\def\Q{{\mathbb Q}}
\def\R{{\mathbb R}}
\def\sss{{\mathbb S}}
\def\T{{\mathbb T}}
\def\Z{{\mathbb Z}}

\def\cA{{\mathcal A}}
\def\cC{{\mathcal C}}
\def\cF{{\mathcal F}}
\def\cO{{\mathcal O}}
\def\cP{{\mathcal P}}
\def\cR{{\mathcal R}}

\def\spzb{{\overline{\Spec\,\Z}}}

\def\HC{\text{HC}}
\def\bfh{{\bf H}}

\def\Mod{{\bf Mod}}

\def\id{{\mbox{Id}}}

\def\hatz{{\hat\Z^\times}}

\def\rmax{{\R_+^{\rm max}}}
\def\zmax{\Z_{\rm max}}

\def\Spec{{\rm Spec\,}}

\def\spm{{\sss[\pm 1]}}

\def\End{{\rm End}}
\def\Hom{{\mbox{Hom}}}

\def\nt{\N^{\times}}

\def\fr{{\rm Fr}}

\def\lbt{{\tilde \Lambda}}
\def\lbto{{\lbt^{\rm o}}}

\def\Ses{{\Se_*}}
\def\Se{\frak{Sets}}
\def\fin{\frak{Fin}}
\def\gop{{\Gamma^{\rm o}}}

\def\totgs{oriented groupo\"ids }
\def\totg{oriented groupo\"id }

\def\frg{\frak{g}}

\begin{document}

\title{Cyclic theory and the pericyclic category}

\author{Alain Connes}
\address{Coll\`ege de France and IHES, France}
\email{alain@connes.org}

\author{Caterina Consani}
\address{Department of Mathematics, Johns Hopkins University,
Baltimore MD 21218, USA}
\email{cconsan1@jhu.edu}
\thanks{The second author is partially supported by the Simons Foundation collaboration grant n. 691493.}

\date{\today and, in revised form, \today.}

\dedicatory{}


\maketitle

\begin{quote}{\it Pour moi le ``paradis originel" pour
l'alg\`ebre topologique, n'est nullement la sempiternelle cat\' egorie $\Delta$ semi-simpliciale, si utile soit-elle, et encore moins celle des espaces topologiques (qui l'une et l'autre s'envoient dans la 2-cat\' egorie des topos, qui en est comme une enveloppe commune), mais bien la cat\' egorie {\bf Cat} des
petites cat\' egories, vue avec un oeil de g\' eom\`etre
par l'ensemble d'intuitions, \' etonnamment riche, provenant des topos. En effet, les topos ayant comme
cat\' egories des faisceaux d'ensembles les $\widehat\cC$, avec $\cC$
dans {\bf Cat}, sont de loin les plus simples des topos
connus, et c'est pour l'avoir senti que j'insiste tant
sur ces topos
``cat\' egoriques" dans SGA  IV.\newline
\centerline{\hspace*{1in}Alexander Grothendieck}}
\end{quote}

\begin{abstract}
We give an historical perspective on the role of the cyclic category $\Lambda$ in the development of  cyclic theory. This involves a continuous interplay between  the extension in ``characteristic one" and in  $\sss$-algebras, of the  traditional development of cyclic theory, and the geometry of the toposes associated with several small categories involved. We clarify the link between various existing presentations of the cyclic and the epicyclic categories and we exemplify the role of the absolute coefficients by presenting the ring of the integers as polynomials in powers of $3$, with coefficients in $\sss[\pm 1]$. Finally, we introduce the pericyclic category  which unifies and refines  two conflicting notions of epicyclic space existing in the literature.
\end{abstract}

\section{Introduction}
\label{catcat}

The introduction  (by the first author of this article) of cyclic cohomology as the analog of de Rham cohomology in noncommutative geometry, was inspired by the development of the quantized calculus \cite{Ober81, CoExt, CoIHES, Co-book} where cyclic cohomology of noncommutative algebras  is the natural receptacle for the Chern character of K-homology cycles. Both the periodicity operator $S: \HC^n(\cA)\to \HC^{n+2}(\cA)$ and the long exact sequence relating cyclic and Hochschild cohomologies of an algebra $\cA$ stem naturally from the analytic theory. 
After  many repetitive, algebraic manipulations  it became clear \cite{CoExt} that laying at the core of the theory there is a  small category $\Lambda$: the cyclic category, encoding combinatorially the rules of computations for an algebra  through a representation of $\Lambda$ in the vector space $\cA^\natural$ (cyclic module) obtained as the direct sum $\cA^\natural:= \bigoplus \cA^{\otimes n}$ of tensor powers of  $\cA$. This result shows, in particular, that cyclic cohomology and homology are two instances of  derived functors in homological algebra. At the conceptual level,  it  determines an embedding of the category of algebras in the larger abelian category of cyclic modules, thus providing a hint for a motivic development of noncommutative geometry.\newline
The cyclic category $\Lambda$  is, by definition, an extension of the category $\Delta$ of finite ordinals. Its classifying space $B\Lambda$  is the infinite projective space $\mathbb P^\infty(\C)$ \cite{CoExt}. This fact points out the similarity of $\Lambda$ with the compact group $U(1)$, as it is well known that $BU(1)=\mathbb P^\infty(\C)$. However, one also ought to keep in mind that while  shifting our focus from a small category  to the homotopy type of its classifying space  we lose some relevant information. For instance, $B\Delta$ is contractible (the category $\Delta$ has a final object), while the restriction to $\Delta$ of the cyclic representation  in $\cA^\natural$ is known to  govern Hochschild homology.\newline A better way to store geometrically the relevant categorical  information is using the associated presheaf topos, as explained in the quote of A. Grothendieck at the beginning of this article. A small category $\cC$ is geometrically represented  by the (Grothendieck) topos $\widehat \cC$ of contravariant functors $\cC\longrightarrow \Se$ to sets. Both  Hochschild and cyclic homologies are described naturally by derived functors on sheaves of abelian groups on the corresponding toposes $\widehat\Delta$ and $\widehat\Lambda$. The additional structure provided by the $\lambda$-operations determines two extensions of $\Delta$ and $\Lambda$, yielding  to the epicyclic refinement $\lbt$ of the cyclic category \cite{bu1, Loday, topos} which is, in essence,  the semi-direct product of $\Lambda$ by an action of the multiplicative monoid $\nt$ of non-zero positive integers. This action  determines a functor $\Mod:\lbt\longrightarrow \nt$, and at the topos level, a geometric morphism from the epicyclic topos  to the topos $\widehat\nt$ \cite{cyctop}. The notion of ``point" in a topos leads, even for toposes of presheaf type, to remarkable new spaces. For example, the presheaf  topos $\widehat\nt$ \cite{CCas} gave the first  example of a Grothendieck topos whose space of points is noncommutative, \ie a quotient of an ordinary space by an equivalence relation, without any Borel section. The relation between the space of points of $\widehat\nt$ and the ad\`ele class space of the rationals $\Q^\times\backslash \A_\Q$ \cite{Cad} has culminated with the discovery of the topos $[0,\infty)\rtimes \nt$ \cite{CCscal1}, and the result that the space of its points  describes the ``$\zeta$-sector" of $\Q^\times\backslash \A_\Q$: \ie the quotient $\Q^\times\backslash \A_\Q/\hatz$ of the ad\`ele class space of $\Q$ by the right action of $\hatz$. The upgrading of  $\widehat\nt$ to the topos $[0,\infty)\rtimes \nt$ was  motivated by the search for suitable structure sheaves on these toposes providing a precise geometric structure. As  explained in Section \ref{sectprojgeom}, the  algebras involved here are of ``characteristic one" and these  structures appear naturally in the study of the cyclic and epicyclic categories. In Section \ref{sectgeomdev} we show how an independent study of these algebras in ``characteristic one" and the discovery of the existence for these structures of an analog of the Frobenius endomorphisms, led us to the definition of the geometries of the Arithmetic and Scaling sites.\newline 
The search for a universal arithmetic then guided us to the study of algebras in the context of Segal's $\Gamma$-rings. This development is reported in Section \ref{sectuniv} and  illustrates the versatility of the ``paradis originel" for the topological algebras mentioned by Grothendieck. Indeed, the fundamental objects of the relevant tensor category of $\sss$-modules (\ie $\Gamma$-sets) are defined as functors  from a small category $\gop$ (a skeleton of the category of finite pointed sets) to sets.  Moreover, the existence of the absolute base $\sss$ is dictated by the development of topological cyclic cohomology within the realm of ring spectra. As shown in \cite{DGM},   Segal's $\Gamma$ spaces provide the right framework for topological Hochschild and cyclic theories. By the Dold-Kan correspondence they give rise to the tools apt to develop homological algebra for $\sss$-algebras. In Section \ref{sectuniv} we explain how this point of view has brought us to the definition of the structure sheaf for the Arakelov compactification $\spzb$ of the spectrum of the integers and the Riemann-Roch theorem for $\spzb$. In Section \ref{integers} we focus on the resulting presentation of the ring $\Z$ of the integers as polynomials with coefficients in the spherical group ring $\sss[C_2]=\spm$, where $C_2$ is the cyclic group of order $2$. This algebra  also appears as global sections of the structure sheaf of $\spzb$.\newline 
Finally,  Section \ref{sectbigperi} introduces the main novelty of the present paper namely the {\em pericyclic category} $\Pi$. This  category unifies and refines  two conflicting notions of epicyclic space existing in the literature. 
 The epicyclic category of Section \ref{sectcat} was originally introduced by T. Goodwillie in a letter to F. Waldhausen and it was implemented in \cite{bu1} to define the corresponding notion of epicyclic space.  Later on, Goodwillie introduced in \cite{GW} a  topological category denoted $S\cR\cF$, also called ``epicyclic", which is at the root of topological cyclic homology. The pericyclic category $\Pi$ (see \S \ref{sectpericyc}) is the combinatorial version of a slightly different version $S'\cR\cF$ of $S\cR\cF$. We show that the two notions of epicyclic space mentioned above are encoded by two covariant functors $\Pi\longrightarrow \Se$ and we expect that the 
pericyclic category $\Pi$ will play the role of a relevant backbone in our future development of a combinatorial version of cyclic homology (hinted at in Section \ref{sectuniv}).

 \section{Cyclic categories and projective geometry}\label{sectprojgeom}
 
 Our original motivation for studying exotic algebraic structures of ``characteristic one" was to achieve a good algebraic understanding of the cyclic and epicyclic categories, and in particular of the intriguing isomorphism between the cyclic category and its opposite $\Lambda^{\rm o}$. It turns out this isomorphism can be  read as  the transposition of linear maps in projective geometry,  in the natural extension of linear algebra within the context of semifields. The conceptual reason behind the algebraic interpretation of the self-duality of 
 the cyclic category  requires some  explanation. Originally,  $\Lambda$ has been defined in terms of homotopy classes of degree one, non-decreasing self-maps of the circle $S^1$, which are compatible with roots of unity  \cite{CoExt}. A better interpretation of it, (see \cite{FT}), is obtained by lifting these maps to the universal cover $\R$ of $S^1$, and  then by restricting these maps  to $\Z$. In this way,  the morphisms between the objects in $\Lambda$ are described by  non-decreasing maps $\phi:\Z\to \Z$ fulfilling the equality, for some $n, m\in \N$ 
\begin{equation}\label{linn}
 \phi(x+n)=\phi(x) + m\qquad \forall x\in \Z.   \end{equation}
 An equivalence between two such maps $\phi, \phi'$ is expressed by the following relation
 \begin{equation}\label{linn1}
 \phi \sim \phi'\iff \exists k\in \Z~\text{s.t.}~ \phi'(x)=\phi(x)+km \qquad  \forall x\in \Z. 
  \end{equation}  
  We shall see that one   may interpret the non-decreasing condition together with \eqref{linn}   as a  linearity property, and  \eqref{linn1} as a projective equivalence. This  is made possible by appealing to the theory of semirings and semifields which is the natural development of the classical theory of rings and fields when one drops the existence of additive inverses. Semirings arise naturally as endomorphisms of any object in a category with a zero object (\ie initial and final) in which the natural maps from direct sums to direct products are isomorphisms. Moving from fields to semifields introduces a few new objects of great interest: the first one is the Boolean semifield $\B= \{0,1\}$,  the unique finite semifield which is not a field. The addition in $\B$ is idempotent ($1+1:=1$) and a semiring contains $\B$ if and only if it is also idempotent (aka of ``characteristic one"). One also sees that there is only one semifield  whose multiplicative group is infinite cyclic: this is the tropical semifield $\zmax$.   More precisely, $\zmax$ is the set $\Z\cup \{-\infty\}$ endowed with the  operation 
  \begin{equation*}\label{addz}
  n\vee m:= \max(n,m)\qquad \forall n,m \in \Z\cup \{-\infty\}
  \end{equation*}
  replacing the classical  addition in $\Z$, while the addition of the integers defines the multiplication in $\zmax$. In particular, the multiplicative group of $\zmax$ is the infinite cyclic group $\Z$.\newline
Next, for each integer $n>0$, one defines a $\zmax$-module $E_n$  whose underlying set is $\Z\cup \{-\infty\}$ with the same $\vee$ operation, and where the action of $\zmax$ on $E_n$ is defined by  
  \begin{equation*}\label{addz2}
 k\bullet m:=m+kn\qquad \forall k\in \zmax, \ m \in E_n.
  \end{equation*}
One has the following interpretation for  a map $\phi: \Z\to \Z$ 
   \begin{equation*}\label{addz1}
 \phi \ \text{is non decreasing}\ \iff  \phi(n\vee m)= \phi(n)\vee \phi(m)\qquad \forall n,m \in \Z.
  \end{equation*}
Therefore,  the non-decreasing maps $\phi:\Z\to \Z$ fulfilling  \eqref{linn} can be seen as the non-zero elements of the set $\Hom_{\zmax}(E_n,E_m)$ which determine,  in view of \eqref{linn1}, projective equivalence classes of ($\zmax$-)linear maps. Thus, one concludes that the cyclic category  $\Lambda$ coincides with the category of projective morphisms between the modules $E_n$, and that the  isomorphism $\Lambda\simeq \Lambda^{\rm o}$ is expressed in terms of  transpositions of such maps \cite{CCproj}. This algebraic description of the cyclic category  hints at the extension from $\Lambda$ to the epicyclic category $\tilde \Lambda$ which is known, by the work of J. L. Loday \cite{Loday}, to govern the $\lambda$-operations  in cyclic homology. In classical projective geometry over fields, morphisms connecting two projective spaces over fields $K$ and $K'$ are constructed from  linear maps $\phi:E\to E'$ of vector spaces, twisted by a morphism $\sigma: K\to K'$ of the fields, that satisfy the rule
\begin{equation*}\label{twist}
\phi(\lambda \xi)=\sigma(\lambda)\phi(\xi)\quad\forall \lambda \in K, ~\xi \in E.
  \end{equation*}
  This suggests the extension of the cyclic category by implementing the endomorphisms of the semifield $\zmax$. In this respect, one of the striking features of an algebra in ``characteristic one" is that for any idempotent semifield $F$ (\ie the addition of any element with itself gives itself)  the map
  \begin{equation}\label{frn}
\fr_a:F\to F, \qquad \fr_a(x):=x^a\quad \forall x\in F
  \end{equation}
  is an injective endomorphism,  for any positive integer $a$ (\cite{G}, Propositions 4.43-4.44). These endomorphisms are  the analogs of the Frobenius endomorphism in finite characteristic, and for $F=\zmax$, they  are the only non-trivial ones.  Using these morphisms to allow for twisted maps between projective spaces amounts to extending \eqref{linn} to 
  \begin{equation*}\label{linna}
 \phi(x+n)=\phi(x) + m\, a\quad\forall x\in \Z. 
  \end{equation*}
  The epicyclic category $\tilde \Lambda$ is described precisely as the extension of the cyclic category  obtained by keeping the same objects $\mathbb P(E_n)$ and implementing new morphisms given by general (twisted) morphisms of projective spaces \cite{CCproj}. Furthermore, there is a natural ``degree" functor $\Mod:\tilde \Lambda\longrightarrow \nt$ which associates to a twisted morphism  by $\fr_a$ its degree $a\in \nt$. The emerging geometric picture   involves geometric morphisms of the toposes dual to the categories.  To promote, when the algebra $\cA$ is commutative, the cyclic module $\cA^\natural$  as a  sheaf of vector spaces over the topos dual to the epicyclic category one proceeds as follows. The forgetful functor which associates to a map of projective spaces $\mathbb P(E_n)\to \mathbb P(E_m)$ the underlying morphism  $\Z/n\Z\to \Z/m\Z$ of sets of points, connects $\tilde \Lambda$ to the category $\fin$ of finite sets. When an algebra $\cA$ over a field $k$ is commutative, the cyclic module $\cA^\natural$ extends to a covariant functor from $\fin$ to vector spaces. To understand this functor as a sheaf of vector spaces of a topos, one considers the topos associated with the small category  $\lbto$ (the opposite category of $\lbt$) which maps to $\nt$ by  the ``degree" functor
 \begin{equation*}
\xymatrix{
 \lbto \ar[d]_(.40){\text{\bf{Mod}}}    \\
\nt 
}
\end{equation*}
This functor extending $\cA^\natural$  governs the $\lambda$-operations in cyclic homology and is discussed in details in \S\ref{sectpericyc}. In \cite{CC6}, these operations are implemented to obtain the cyclic cohomological interpretation of Serre's local $L$-factors  of arithmetic varieties.

\section{Geometric developments}\label{sectgeomdev}

 The main role played by  $\zmax$  in the understanding of the cyclic and epicyclic categories, and  the canonical action \eqref{frn} of $\nt$ by Frobenius endomorphisms on this semifield, led us to endow the ``base"  topos $\widehat \nt$ with the structure sheaf prescribed exactly by this  action. The outcoming geometric structure is the {\em Arithmetic Site} \cite{CCas} which, despite being defined using only countable data of combinatorial nature, admits a one-parameter family of correspondences whose composition law reflects a fundamental action of the multiplicative group $\R_+^*$ of positive real numbers. This  is the group of Galois automorphisms of the semifield $\rmax$ of positive real numbers under the two operations of max for addition and usual multiplication, which plays a central role both in the semiclassical limit of quantum mechanics \ie in idempotent analysis, and tropical geometry. The space of points of the Arithmetic Site over $\rmax$ is canonically isomorphic to the $\zeta$-sector of the ad\`ele class space $\Q^\times\backslash \A_\Q$ of the rationals \cite{CCas}, and  the Galois action of  $\R_+^*$ on this space reflects geometrically the action of the id\`ele class group on the $\zeta$-sector, in complete analogy with the geometric development in finite characteristic. The extension of scalars to $\rmax$ defines a new geometric object: the  {\em Scaling Site} \cite{CCscal1}. As a topos, this is   the semi-direct product $[0,\infty)\rtimes \nt$, where $[0,\infty)$ is the euclidean half line $[0,\infty)$. The points of the topos  $[0,\infty)\rtimes \nt$  coincide with the  ``$\zeta$-sector" already selected by the points of the Arithmetic Site over the same scalar extension. The new additional structure inherited by the Scaling Site is the structure sheaf made by piecewise affine, convex functions. A Riemann-Roch theorem analogous to that holding for elliptic curves is verified for the restriction of the Scaling Site to its periodic orbits  \cite{CCscal1}. A new striking feature is that the dimensions of the cohomologies involved are real numbers and they correspond, in the world of ``characteristic one", to the real dimensions of Murray and von Neumann.  This Riemann-Roch formula  implements only $H^0$ and $H^1$ for divisors. One uses Serre duality to define $H^1$, by-passing in this way, the development of a full cohomological machine. The module $H^0$ of global sections is naturally filtered using the $p$-adic valuation for the periodic orbit associated with the prime $p$, and  this filtration gives rise to  continuous dimensions.\newline The development in \cite{CCscal3} of homological algebra in ``characteristic one"  has shown that the lack of  additive inverses creates considerable difficulties in the process of passing from abelian categories to their counterpart, for idempotent semirings. 
 
  \section{Universal cyclic theory}\label{sectuniv}
 
  Surprisingly, one can circumvent the open problem of the elaboration of homological algebra for idempotent structures by developing cyclic homology in a world parallel to that of noncommutative geometry! This is the world  of algebraic $K$-theory as promoted by D. Quillen and F. Waldhausen. We see it parallel to noncommutative geometry since it is transverse to the standard classification of mathematics into distinct branches. Cyclic homology was dragged into this context at the beginning of the $80$'s, and then developed for spectra,  in the world of ``brave new rings". At the Oberwolfach meeting in $1987$, which was organized to sponsor interactions between functional analysis and algebraic topology, T. Goodwillie \cite{Ober87} emphasized that in the world of ``brave new rings", the ring $\Z$ of the integers becomes an algebra over a more fundamental ring,  namely the sphere spectrum $\sss$.  While the homotopical constructions  are not particularly appealing to   functional analysts,  as iterated loop spaces may get quickly out of explicit control, a more concrete and basic version of the theory has  gradually emerged and it is accounted for in  \cite{DGM}. The  framework in question is that of G. Segal $\Gamma$-spaces which provides a workable model of connective spectra. \newline Our independent viewpoint on this subject is that the theory of  $\Gamma$-spaces  supplies the required  constructions to produce the basic homological algebra, available through the Dold-Kan correspondence,  for a very concrete algebraic set-up solely based on the Segal category $\gop$. This category is a skeleton of the category $\fin_*$ of pointed finite sets. Its main role is to provide the broadest  generalization of a commutative and associative operation. These two properties combine to give a meaning to a ``finite sum" and to obtain a covariant functor $\cF:\fin_*\longrightarrow \Se_*$, where the  image of a map $f:X\to Y$ of finite pointed sets uses the finite sum over the inverse image $f^{-1}(y)$ to obtain the value at $y\in Y$.  This construction  generalizes broadly the Eilenberg-MacLane functor which replaces an abelian group A with the covariant functor $\text{HA}: \fin_*\longrightarrow \Se_*$ that assigns to a finite pointed set $X$ the pointed set of A-valued divisors on $X$ (taking the value $0$ at the base point). Divisors push forward by the formula
   \begin{equation}\label{divisors}
 f_*(D)(y):=\sum_{f(x)=y} D(x).
  \end{equation}  
  The Eilenberg-MacLane functor determines a fully faithful inclusion  of abelian groups into the category of (covariant) functors $\gop\longrightarrow \Se_*$, where morphisms are natural transformations of functors. Then, by implementing an idea of G.~Lydakis, one defines the smash product of such  functors starting from the closed structure determined by the smash product in $\fin_*$. The identity functor $\id: \gop\longrightarrow \Ses$ defines  the easiest example of such an algebra  in this tensor category and it describes  the sphere spectrum $\sss$ in its most elementary form. This is the fundamental brave new ring advocated by Goodwillie  in Oberwolfach.  \newline
  In \cite{CCprel} we show that this theory of $\sss$-algebras embodies all the previous attempts to extend ordinary algebra to reach a  ``universal arithmetic".\newline 
  The main difficulty that one encounters in trying to extend the structure sheaf of $\Spec \Z$ by taking into account the point at $\infty$ (the archimedean place), comes from the lack of an analog, for the archimedean prime, of the subring of $\Q$ which, for a non-archimedean prime $p$, is  the localization of $\Z$ at $p$. In $p$-adic terms this corresponds to elements of $\Q$ which belong to the local ring $\Z_p$ of $p$-adic integers (whose  $p$-adic norm is $\leq 1$). At the archimedean place, one meets the difficulty that the condition $\vert q\vert \leq 1$ does not define a subgroup of $\Q$. However, this problem is successfully  solved by the theory of $\sss$-algebras. Indeed,  the Eilenberg-MacLane functor H replaces faithfully a ring R by the $\sss$-algebra HR and this applies in particular to the rings $\Z$ and $\Q$. In this categorical universe the $\sss$-algebra H$\Q$ admits the non-trivial subalgebra $\cO_\infty: \gop\longrightarrow \Ses$ defined by the condition on $\Q$-valued divisors:
  \begin{equation}\label{divisorsq}
 \sum_{x\in X}\vert  D(x)\vert \leq 1.
  \end{equation} 
   This condition is  stable by push forward  and as such it defines an $\sss$-module. Moreover, since the product map 
    \begin{equation*}\label{product}
\cO_\infty(X)\wedge \cO_\infty(Y)\to \cO_\infty(X\wedge Y)
  \end{equation*} 
  is the restriction of the product map for H$\Q$ which respects the basic condition \eqref{divisorsq}, one sees that $\cO_\infty$ is  an $\sss$-subalgebra of H$\Q$. This fact allows one  to define the structure sheaf of $\spzb$ as a subsheaf of the constant sheaf $\Q$ \cite{CCprel}.\newline 
  A relevant application of this construction is the detection of the Gromov norm, which is obtained through a refinement of the ordinary homology with rational coefficients using $\cO_\infty$  as $\sss$-module \cite{CCgromov}. \newline 
  With the point at infinity taken into account, $\spzb$ behaves like a projective curve, and the global sections of its structure sheaf form a very simple extension of the basic algebra $\sss$. This is the spherical group ring $\sss[C_2]=\spm$, which associates to a finite pointed set $X$ the smash product of $X$ with the pointed, multiplicative monoid $\{-1,0,1\}$. By working over this base, one gains the categorical interpretation of  additive inverses! For an abelian group A, for example, the $\sss$-module HA is naturally an $\spm$-module using the rule $(-1)\times x=-x$ for all $x\in \text{A}$.\newline
  This construction seems to have the potential to go quite far. In the recent article \cite{RRinv} for example, we show how to parallel, on $\spzb$, the adelic proof of the Riemann-Roch formula for function fields given by A. Weil, including the use of Pontryagin duality.  The  definition of the  cohomologies is based on the universal arithmetic explained above and on an extension for $\sss$-modules of the Dold-Kan correspondence. To an Arakelov divisor $D$ on $\spzb$ one associates a $\Gamma$-space $\bfh(D)$ using a short complex of $\spm$-modules directly related to the divisor. The two cohomologies $H^0(D)$ and $H^1(D)$ are then defined {\it independently}, as well as the notion of their dimension over $\spm$. The Riemann-Roch formula also proven in \opcit  equates the {\it integer valued} Euler characteristic of a divisor $D$ with a simple modification of the traditional expression (\ie the degree of the divisor plus log 2). The integer-valued topological side of the formula involves, besides  the ceiling function,  the division by log 3.  This result unveils the special role of the number $3$ among the integers and  the triadic expansion of the integers as polynomials $P(X)$,  with coefficients in $\{-1,0,1\}$ evaluated at $X=3$. In Section \ref{integers} we shall explore the extent to which one can view  $\Z$ as a polynomial ring over $\spm$. \newline A long-range goal is the computation of  the universal  cyclic cohomology for $\spzb$.

 \section{$\spm[3]$ and the ring of integers}\label{integers}
 
 In our recent work \cite{RRinv} on the arithmetic Riemann-Roch theorem  for $\spzb$, the number $3$ plays the  role of a natural generator, allowing one to label the integers $m\in \Z$ as polynomials in powers of $3$ with coefficients in $\spm$
 \begin{equation}\label{polynomial}
 P(X)=\sum_{j=0}^k a_j X^j, \ \ a_j\in  \{-1,0,1\},	\ \forall j.
 \end{equation}
 The subtlety of the description of the ring structure of the integers $\Z$ in this parametrization derives from    the addition of  polynomials with coefficients in the multiplicative monoid $\{-1,0,1\}$. Here $\{-1,0,1\}$ is endowed with the obvious multiplication rule, thus one knows how to multiply monomials \ie $a\,X^n\times b\, X^m=(ab)\,X^{n+m}$, while  the product of polynomials is obtained  uniquely, using the distributive law, provided one knows how to add them. The sum of two monomials $a\,X^n , b\, X^m$ of different degrees is simply the polynomial $a\,X^n + b\, X^m$ and when the degrees are the same the distributivity law $a\,X^n + b\, X^n=(a+b)X^n$ reduces the reconstruction of the sum to polynomials of degree $0$. \newline
 We shall make sure that addition is defined on polynomials of degree $0$ in such a way that the following rule holds:\vspace{.03in}
 
 \centerline{\bf The sum of $3$ polynomials of degree $0$ is a polynomial of degree $\leq 1$.}\vspace{.03in}
 
\noindent  Using this hypothesis it follows, by induction on the degree $n$, that the sum of two polynomials of degree $\leq n$ is a polynomial of degree $\leq n+1$. The need to require that the sum of three polynomials of degree $0$ is of degree $\leq 1$ (rather than for just the sum of two), is due to the role of the  carry over in adding arbitrary polynomials. \newline
  Let us now consider the addition of two elements of $\{-1,0,1\}$. If one of them is $0$ or if they have opposite signs, then the sum is the obvious one. Thus,   using the symmetry $x\mapsto -x$, we can assume that they are both  equal to $1$.  
 The only new rule that we add is the following 
 \begin{equation}\label{addition} 
 1+1 = X-1.	
 \end{equation}
 By applying \eqref{addition}, one   adds two polynomials of degree $0$ using the following table
\[
\begin{tabular}{c| c c c}
    + & $-1$ & 0 & 1 \\\hline\rule{0pt}{3ex}
    $-1$ & $1-X$ & $-1$ & 0 \\
    0 & $-1$ & 0 & 1 \\
    1 & 0 & 1 & $X-1$ \\
\end{tabular}
\]
To specify the sum of three polynomials of degree $0$, it is enough to know how to add to the polynomial $X-1$  the elements of $\{-1,0,1\}$. Using the associativity of the sum and \eqref{addition}, one has
\[
(X-1)+1=X, \qquad  (X-1)-1=X-(1+1)=X-(X-1)=1.
\]
It follows then that, as required,  the sum of  $3$ polynomials of degree $0$ is a polynomial of degree $\leq 1$. 
 The following result determines the sum of any pair of polynomials under \eqref{addition}
 
\begin{proposition}\label{addpol} Let $P(X)$ and $Q(X)$ be two polynomials as in \eqref{polynomial}, of degree at most $n$. Then there exists a unique polynomial of degree at most $n+1$ which coincides, using the rule \eqref{addition}, with the sum $P(X)+Q(X)$.
\end{proposition}
\proof Let assume (by induction) that the statement holds for any two polynomials of degree $<n$ and let consider a pair $P(X)$, $Q(X)$ both of degree $n$. We decompose $P(X)=P_1(X)+\alpha_n X^n$,  $Q(X)=Q_1(X)+\beta_n X^n$, with $\alpha_n,\beta_n\in\{-1,0,1\}$. Then
\[
P(X)+Q(X)=P_1(X)+Q_1(X)+(\alpha_n+\beta_n)X^n,
\]
where $P_1(X)+Q_1(X)$ is, by induction, the sum of a polynomial of degree $n-1$ with a term $\gamma X^n$, $\gamma\in\{-1,0,1\}$. Then, from the specification of the sum of three polynomials of degree zero we have $\alpha_n+\beta_n+\gamma=\delta+\epsilon X$, with $\delta,\epsilon\in \{-1,0,1\}$.  Thus the sum $P(X)+Q(X)$ can be written as a polynomial of degree at most $n+1$, where the terms of degree $\geq n$ are $\delta X^n+\epsilon X^{n+1}$.\endproof 

We obtain the following intriguing conclusion

\begin{proposition}\label{ringpol} The set of polynomials as in \eqref{polynomial}, under the addition stated in Proposition \ref{addpol}, and the unique associated product, forms a ring isomorphic to $\Z$. 
\end{proposition} 
\proof The map 
\[
\theta: \sum \alpha_j X^j\mapsto \sum \alpha_j 3^j
\]
is a bijection from the set of polynomials as in \eqref{polynomial}  with the integers. This map  is compatible with addition and multiplication in view of Proposition~\ref{addpol}. \endproof 

There is a function that allows one to compute effectively sums of polynomials as above. This is a function  of three variables $\alpha,\beta,\gamma$  in $\{-1,0,1\}$ which specifies their sum   as a polynomial of degree at most $1$ in $X$. To determine  $\alpha+\beta+\gamma$ one can  think of it as a sum of $3$ real numbers and write the set $\{-1,0,1\}+\gamma$ where this sum belongs, to obtain the coefficient of $X$. At the algebraic level,   $\alpha+\beta+\gamma$ can be written explicitly using the finite field $\F_3=\{-1,0,1\}$ as follows 
\[
\alpha+\beta+\gamma=\overline{\alpha+\beta+\gamma}+ 
\overline{\sigma(\alpha,\beta,\gamma)}\ X
\]
where  $\overline m$ is the residue of $m$ modulo $3$ and 
\[
\sigma(\alpha,\beta,\gamma):=\alpha \beta \gamma-\alpha^2 \beta-\alpha^2\gamma-\alpha \beta^2-\alpha \gamma^2-\beta^2 \gamma-\beta \gamma^2.
\]	
This determines  the algebraic expression for the sum of three polynomials of degree $0$ as in the proof of Proposition~\ref{addpol}. Then, by induction on  $n$, one  obtains an algebraic expression for the coefficients of $X^n$ in    the sum $\sum_j\alpha_jX^j+\sum_j\beta_jX^j$ of two polynomials. Indeed, one defines by induction  
 polynomials \sloppy $s_n(\alpha_0,\ldots,\alpha_{n-1},\beta_0,\ldots,\beta_{n-1})$ with coefficients in $\F_3$ as follows
\[
s_0:=0,\quad s_n=\sigma(\alpha_{n-1},\beta_{n-1},s_{n-1}).
\]
Then the coefficient of $X^n$ in the sum $\sum_j\alpha_jX^j+\sum_j\beta_jX^j=\sum_j\gamma_jX^j$ is given by 
\begin{equation}\label{witthold}
	\gamma_n=\alpha_n+\beta_n+s_n.
\end{equation}
This construction provides a sequence of polynomials $\gamma_n(\alpha_1,\ldots,\alpha_n,\beta_1,\ldots,\beta_n)$ and finally one has 
\begin{lemma}
The  polynomial $\gamma_n(\alpha_1,\ldots,\alpha_n,\beta_1,\ldots,\beta_n)$ vanishes provided the  components  $\alpha_n,\alpha_{n-1}$ and $\beta_n,\beta_{n-1}$ all vanish.	
\end{lemma}
\proof It is enough to show that $s_n=0$ if $\alpha_{n-1}=0$ and $\beta_{n-1}=0$. This follows from the inductive definition of $s_n=\sigma(\alpha_{n-1},\beta_{n-1},s_{n-1})$. \endproof 

Since one works over $\F_3$  $x^n=x$ if $n$ is an odd number, while $x^n=x^2$ if $n$ is a non-zero even number. Thus one can reduce the  exponents of the variables in the polynomials $\gamma_n$ to be  at most $2$. 
One has $s_1(\alpha_0,\beta_0)=-\alpha_0 \beta_0^2-\alpha_0^2 \beta_0$, while $s_2$ in reduced form is given by
\begin{multline*}
s_2(\alpha_0,\alpha_1,\beta_0,\beta_1)=\alpha_1 \alpha_0^2 \beta_0^2+\alpha_0^2 \beta_0 \beta_1^2+\alpha_1^2 \alpha_0^2 \beta_0+\alpha_0^2 \beta_0^2 \beta_1-\alpha_1 \alpha_0^2 \beta_0 \beta_1\\+\alpha_1^2 \alpha_0 \beta_0^2+\alpha_0 \beta_0^2 \beta_1^2+\alpha_1 \alpha_0 \beta_0-\alpha_1 \alpha_0 \beta_0^2 \beta_1+\alpha_0 \beta_0 \beta_1-\alpha_1 \beta_1^2-\alpha_1^2 \beta_1
\end{multline*}
To write in full $s_3(\alpha_0,\alpha_1,\alpha_2,\beta_0,\beta_1,\beta_2)$ would take several lines.\vspace{.05in}

Conceptually, one  sees that the above construction is described by the addition rule of two Witt vectors over $\F_3$. The number $3$ is the only prime for which the Witt vectors with only finitely many non-zero components form an additive subgroup of the Witt ring. More precisely one can state the following 

\begin{proposition}\label{ringpolbis} Witt vectors with only finitely many non-zero components form a subring of the ring  $\mathbf W(\F_3)$. This subring is isomorphic to $\Z\subset \Z_3$. 
\end{proposition}
\proof  For any prime $p$ the  map 
\[
\mathbf W(\F_p)\ni \xi=(\xi_j)\mapsto \tilde \tau(\xi):= \sum_{j=0}^{\infty} \tau(\xi_j) p^j,
\]
where $\tau$ is the Teichm\"uller lift, is known to be a ring isomorphism between the Witt ring $\mathbf W(\F_p)$ and the ring $\Z_p$ of $p$-adic integers. For $p=3$,   $\tau(\F_3)=\{-1,0,1\}\subset \Z\subset \Z_3$, thus $\tilde \tau$ maps Witt vectors with only finitely many non-zero components into the subring $\Z\subset \Z_3$. This map is an injective ring homomorphism and it surjects onto $\Z$ since by Proposition~\ref{ringpol}  any integer can be written as a finite sum of powers of $3$ with coefficients in $\{-1,0,1\}$. \endproof 
 
 \section{The cyclic and epicyclic categories as subcategories of {\bf Cat}}\label{sectcat}
  
  In Section \ref{sectprojgeom} we reviewed the cyclic and the epicyclic categories using  concepts of linear algebra in characteristic one, in particular we have described the objects of these categories as projective spaces over the semifield $\zmax$. Given these results it seems natural to wonder about the existence of further presentations of the same categories: in particular whether their objects and morphisms may be realized within the ``paradis originel" provided by the category {\bf Cat} of small categories. This idea was developed by D. Kaledin in \cite{Kaledin} and it turns out to be a slight variation of our earlier viewpoint developed in \cite{topos}. 
   In that paper we described the cyclic and epicyclic categories in terms of oriented groupoids: this result is reviewed here below. This section is dedicated to comparing our viewpoint with that of \cite{Kaledin} that has the advantage to interpret the collection of the objects as  small categories (and as such objects of {\bf Cat}) and  the morphisms as functors, thus presenting the two categories as subcategories of {\bf Cat}.
   
  \subsection{Oriented groupoids}\label{sectorgroupoids}

We recall that a groupo\"id $G$ can be seen as a small category in which  every morphism is invertible. 
We denote by $G^{(0)}$  the set of objects of $G$ and  by
$r,s:G\to G^{(0)}$ {\it resp.} the range and the source of the morphisms of $G$. We view $G^{(0)}\subset G$ as the subset  of the identity morphisms of $G$. 

The following  definition is a direct generalization of the notion of the right ordered group for a groupo\"id $G$. Given $X\subset G$, we let $X^{-1}:=\{\gamma^{-1}\mid \gamma\in X\}$. 

\begin{definition}  An \totg  $(G,G_+)$ is a groupo\"id $G$ endowed with a subcategory $G_+\subset G$, fulfilling the following relations
\begin{equation}\label{orddefn}
G_+\cap G_+^{-1}=G^{(0)}, \qquad G_+\cup G_+^{-1}=G.
\end{equation}
\end{definition}

Let $H$ be a group acting on a set $X$. The semi-direct product $G:=X\ltimes H$ is a groupo\"id with source, range and composition law defined respectively as follows
\[
s(x,h):= x, \quad r(x,h):= hx, \quad (x,h)\circ (y,k):=(y,hk). 
\]
As it happens in any groupo\"id, the composite $\gamma\circ \gamma'$ of two elements of $G$ is only defined when $s(\gamma)=r(\gamma')$: here this holds if and only if  $x=ky$.\newline
There is a canonical homomorphism of groupo\"ids $\rho:G\to H$, $\rho(x,h)=h$. We recall from  \cite{topos} (Lemma~2.2) the following

\begin{lemma} \label{examsemi0} Let $(H,H_+)$ be a right ordered group. Assume that  $H$ acts on a set $X$. Then the semi-direct product $G=X\ltimes H$, with $G_+:=\rho^{-1}(H_+)$, is an \totg\!\!.
\end{lemma}

We shall consider, in particular, the right ordered group $(H,H_+)=(\Z,\Z_+)$ acting by translation on the set
$X=\Z/(m+1)\Z$ of integers modulo $m+1$. Then, by applying the above lemma one obtains the \totg 
\begin{equation*}\label{examsemi}
\frg(m):=\big(\Z/(m+1)\Z\big)\ltimes \Z.
\end{equation*}

The relation with the epicyclic category is provided by the following statement (see \cite{topos} Corollary 3.11)

\begin{proposition}\label{epicor} The epicyclic category $\lbt$  is canonically isomorphic to the category whose objects are the \totgs $\frg(m)$, $m\geq 0$, and the morphisms are the non-trivial morphisms of \totgs\!\!. 

\noindent The functor which  associates to a morphism of \totgs its class up-to equivalence (of small categories) corresponds, under the canonical isomorphism, with the functor $\Mod:\lbt\longrightarrow \nt$ sending a semilinear map of semimodules over $\F=\Z_{\rm max}$ to the corresponding injective endomorphism $\fr_n\in \End(\F)$ \cf ~\cite{CCproj}.
\end{proposition}

\subsection{Comparison with  \cite{Kaledin}}\label{sectkaledin}
We recall the presentation of the cyclic category given in \cite{Kaledin}: the following notations and the three definitions are taken literally  from \opcit

Let $[1]_\Lambda$ be the category with one object and whose endomorphisms are determined by non-negative integers $a\in \Z_+$. The inclusion $\Z_+\subset \Z$ induces a functor $[1]_\Lambda\longrightarrow {\bf pt}_\Z$. By the Grothendieck construction, sets with a $\Z$-action  correspond to small categories discretely bi-fibered over ${\bf pt}_\Z$, thus  a $\Z$-set $S$ induces, by pullback,  a discrete bi-fibration $[1]_\Lambda^S\to [1]_\Lambda$. For any integer $n\geq 1$, let $[n]_\Lambda$ be the category corresponding to the $\Z$-set $\Z/n\Z$. The category $[n]_\Lambda$ is equivalent to the additive monoid $\N$ of non-negative integers. 
The degree of a functor is defined in terms of the corresponding morphism of monoids
\begin{definition}\label{kale1} A functor $f:[n]_\Lambda\longrightarrow [m]_\Lambda$ is said to be vertical if it is a discrete bi-fibration,  horizontal if $\deg(f)=1$ and non-degenerate if $\deg(f)\neq 0$.
	
\end{definition}
\begin{definition}\label{kale2} The {\em cyclotomic category} $\Lambda R$ is the small category whose objects $[n]$ are indexed by positive integers and whose morphisms $[n]\to [m]$ are given by non-degenerate functors $f:[n]_\Lambda\longrightarrow [m]_\Lambda$.
\end{definition}

\begin{definition}\label{kale3} The cyclic category  $\Lambda$ is the subcategory of $\Lambda R$ with the same set of objects  and whose  morphisms are the horizontal ones. 	\end{definition}

Next, we  relate the above definitions (taken literally from \cite{Kaledin}) the description of the epicyclic category given in \cite{topos} and here recalled in \S \ref{sectorgroupoids}.  
There is a natural functor $\cP$  which associates to an oriented groupo\"id $G$ the  subcategory $G_+=\cP(G)$, viewed as a small category. For $G=\frg(n)$, one has 
\begin{equation}\label{kalcat}
\cP(\frg(n))=[n+1]_\Lambda \qquad \forall n\geq 0.
\end{equation}
The link with the epicyclic category is  then provided by the following

\begin{proposition}
	The functor $\cP$ identifies the epicyclic category $\lbt$ (as described in Proposition \ref{epicor}) with the cyclotomic category $\Lambda R$.  From the commutativity of the diagram
\[
\xymatrix{
    \lbt \ar[r]^{\cP} \ar[d]_\Mod  & \Lambda R \ar[d]^\deg \\
    \nt \ar[r]_{\rm{Id}}& \nt
  }
  \]
the degree as in \cite{Kaledin} corresponds to the module as in Proposition \ref{epicor}. 
\end{proposition}
Since the cyclic category $\Lambda$ is the subcategory of  $\lbt$ with the same objects, and morphisms  fulfilling $\Mod(f)=1$, it follows from the above identification that Definition~\ref{kale3} determines the cyclic category $\Lambda$ inside $\lbt$.

\section{The pericyclic category $\Pi$}\label{sectbigperi}

In \S \ref{sectkaledin} it has been shown that the epicyclic category $\lbt$ can be identified with the subcategory $\Lambda R$ of {\bf Cat} (the name cyclotomic category for $\Lambda R$ seems, therefore, to be redundant and possibly also creating confusion).
Another source of confusion is that the  name {\em epicyclic} is used in the literature with two different meanings. The epicyclic category $\lbt$ of Section \ref{sectcat} was originally introduced by T. Goodwillie in a letter to F. Waldhausen, however later on  Goodwillie introduced in \cite{GW} another category, also called ``epicyclic" (denoted with $S\cR\cF$  \cf~ \cite{DGM}) and here  described in \S \ref{sectgoodw}. In the same \S \ref{sectgoodw} we also correct a ``typo", sometimes reported in the literature, regarding the composition rule for morphisms in the category $S\cR\cF$. In \S \ref{sectspectra} we show that using the category $S'\cR\cF$  (with the modified rule implemented),  a cyclotomic spectrum gives rise to a functor from $S'\cR\cF$ to spectra. In  \S \ref{sectpericyc} we introduce the {\em pericyclic category} $\Pi$, as the combinatorial version of $S'\cR\cF$ and prove that the two notions of epicyclic space (in relation to the two epicyclic categories mentioned above) are encoded in terms of (covariant) functors $\Pi\longrightarrow \Se$. Finally, \S \ref{sectR}  investigates the points of the topos dual to a small category $\cR$ constituent of the pericyclic category. 

\subsection{The  category $S\cR\cF$}\label{sectgoodw}

The small category $\cR\cF$ (see \cite{DGM} Definition 6.2.3.2)  has one object $a$ for each integer $a>0$ and for any pair of objects $a,b$ the morphisms connecting them are
\[
\Hom_{\cR\cF}(a,b):=\{f_{r,s}\mid r,s\in \nt, a=rbs\}.
\]
In particular one has:  $\Hom_{\cR\cF}(a,b)=\emptyset$ unless $b\vert a$. The composition of morphisms is assigned by
\[
f_{r,s}\circ f_{p,q}=f_{rp,sq}.
\]
Thus one sees that the following maps define functors  sending all the objects to the unique object $\bullet$ of $\nt$ and acting on morphisms as follows

\begin{equation}\label{resfrob}
{\rm Res}:\cR\cF\longrightarrow \nt, \  \ {\rm Res}(f_{p,q}):=p\in \nt, \ \  \ \fr:\cR\cF\longrightarrow \nt, \  \ \fr(f_{p,q}):=q\in \nt.
\end{equation}

The small category $S\cR\cF$ (see \opcit Remark 6.2.3.5) is topological and has the same objects as $\cR\cF$ but  implements an angular parameter $\theta$ 
in the morphisms (of $\cR\cF$),   viewed as an element of $\T=\R/\Z$ or equivalently understood as the element $e(\theta):=\exp(2\pi i \theta)$ of $U(1):=\{z\in \C\mid \vert z\vert =1\}$. An {\em epicyclic space} is a functor from  $S\cR\cF$ to topological spaces. In this framework the first index $r$ in the morphisms $f_{r,s}$ corresponds to  restriction maps (see \cite{DGM}, 6.2.3) while the second index labels the {\em inclusions of fixed points}.\newline
 Let us carefully review the role of such inclusions. Given a $\T$-space $X$, \ie a space $X$ endowed with an action $\theta\mapsto \rho_\theta\in {\rm Aut}(X)$ of $\T$ on $X$, one associates to an integer $q>0$ the cyclic subgroup $C_q=\{\theta \in  \T\mid q\theta =0\}$ (equivalently the group of $q$-th roots of unity in $U(1)$) and the subspace  $X^{C_q}$ of the $C_q$ fixed points. The  canonical inclusion $i_q:X^{C_q}\to X$  is $\T$-equivariant by construction, however, because of the canonical isomorphism 
\begin{equation*}\label{good1}
\alpha: \T\to \T/C_q, \quad \alpha(\theta)=\frac 1q \theta, \quad \alpha^{-1}:\T/C_q\to \T,\quad \alpha^{-1}(\theta)=q\theta 	
\end{equation*}
one modifies the action of $\T$ on $X^{C_q}$ by composing  with $\alpha$. For this new action $\rho'=\rho\circ \alpha$ on $X^{C_q}$ fulfilling $\rho'_{q\theta}=\rho\circ \alpha(q\theta)=\rho_{\theta}$, one has
\begin{equation}\label{good2}
	\rho_{\theta}i_q(x)=i_q(\rho'_{q\theta}(x)).
\end{equation}
Next, we review in detail the notion of an epicyclic space given in \cite{DGM} (Definition 6.2.3.1).
For any integer $k>0$, the edgewise subdivision is an endofunctor ${\rm sd}^k:\Delta\longrightarrow \Delta$ obtained by concatenating $k$ copies of an ordinal: it multiplies the cardinality by $k$. Given a simplicial set $Y$ one lets ${\rm sd}_k(Y):=Y\circ {\text{sd}}^k$ be the simplicial set obtained by pre-composition with these endofunctors. If moreover $Y$ is a cyclic set then the cyclic structure extends to ${\text{sd}}_k(Y)$ provided one replaces the cyclic category $\Lambda$ with the $k$-cyclic  category $\Lambda_k$ defined by replacing the cyclic group $C_n$ involved as automorphisms of the object of cardinality $n$, by the  group $C_{nk}$ with generator $t_{(n,k)}$, while keeping the same relations of the generator with faces and degeneracies. The fixed points $Z^{C_k}$  in a $k$-cyclic set $Z$ form a cyclic set. Thus, given a cyclic set $Y$ and an integer $k>0$ the fixed points ${\text{sd}}_k(Y)^{C_k}$ in the edgewise subdivision ${\text{sd}}_k(Y)$  form a cyclic set.
\begin{definition}\label{defnepicycspace}
	An {\em epicyclic space} $(Y,\phi)$ is  a pointed cyclic space  $Y$ equipped with pointed cyclic maps  $\phi_q:{\rm sd}_q(Y)^{C_q}\to Y$, $\phi_1=\id$, so that the following diagrams commute 
\begin{equation*}\label{epic}
\xymatrix{
{\text{sd}}_k({\text{sd}}_r(Y)^{C_r})^{C_k} \ar[d]_{{\text{sd}}_k(\phi_r)}  \ar[r]^-{\text id} &\ar[d]^{\phi_{kr}}  {\rm sd}_{kr}(Y)^{C_{kr}} \\
{\text{sd}}_k(Y)^{C_k} \ar[r]^-{\phi_k} & Y
}
\end{equation*}
\end{definition}
 
 The geometric realization $\vert Z\vert$ of a $k$-cyclic set $Z$ admits a canonical action of the group $\R/k\Z$, and the fixed points of this action under the subgroup 
$\Z/k\Z\sim C_k$ give rise to the geometric realization of the cyclic set $Z^{C_k}$. Moreover, one has a canonical homeomorphism of geometric realizations (\cite{DGM} Lemmas 6.2.2.1, 6.2.2.2)

\begin{equation*}\label{geomreal}
	D_a: \vert{\text{sd}}_a(Y)\vert\to \vert Y\vert, \qquad D_a(s+a\Z)(y)=\left(\frac sa +\Z\right)D_a(y).
\end{equation*}

Given an epicyclic space $(Y,\phi)$, one introduces for each integer $a>0$ the $\T$-space 
\[Y(a):=\vert{\text{sd}}_a(Y)^{C_a}\vert\simeq \vert Y\vert^{C_{a}}
\]
and the maps $i_q:Y(qa)\to Y(a)$ corresponding to the inclusion of the fixed points 
\[
Y(qa)\simeq \vert Y\vert^{C_{qa}}\subset \vert Y\vert^{C_{a}}\simeq Y(a).
\] 
Then it follows from \eqref{good2} that the compatibility of the maps $i_q$ with the $\T$-action is expressed by the following formula
\begin{equation}\label{good2bis}
	\rho_{\theta}i_q(x)=i_q(\rho_{q\theta}(x)).
\end{equation}
Now we come to the point of issue. The composition of the morphisms is given in \cite{GW} by the formula
\begin{equation}\label{goodw1}
(\theta,f_{r,s})\circ (\tau, f_{p,q})=(\theta+s \tau,f_{rp,sq}).
\end{equation}
However, one should associate to the pair $(\theta,f_{r,s})$ the morphism $\phi_ri_s\rho_\theta$, instead of $\phi_r\rho_\theta i_s$, the reason being that the latter labeling is redundant since $\rho_\theta$ applied to the range of $i_s$ is unchanged when $\theta$ is multiplied by $s$-torsion. Then, standing \eqref{good2bis},   \eqref{goodw1} ought to be replaced with
\begin{equation}\label{goodw1bis}
 (\theta,f_{r,s})\circ (\tau, f_{p,q})=(q\theta+\tau,f_{rp,sq})
\end{equation}
which corresponds to 
\[
i_s\rho_{\theta}i_q\rho_{\tau}=i_si_q\rho_{q\theta}\rho_{\tau}.
\]
Next, we rely on \eqref{goodw1bis} to define the variant $S'\cR\cF$ of the category $S\cR\cF$.
Thus by construction $S'\cR\cF$ is an extension of $\cR\cF$ by $U(1)$, and there is a natural forgetful functor $\sigma:S'\cR\cF\longrightarrow \cR\cF$ which is the identity on objects and maps $(z,f_{r,s})\mapsto f_{r,s}$. We now compare this category with the monoid $U(1)\ltimes \nt$ considered in \cite{bu1}, whose elements are pairs $(z,s)\in U(1)\ltimes \nt$, and the composition law is defined by
\[
(z_1,s)\circ (z_2, q)=(z_1z_2^s,sq).
\]
One obtains the following commutative diagram 
\begin{equation}\label{1o}
\xymatrix{
(U(1)\ltimes \nt)^{o} \ar[d]_p & \ar[l]_-{\text{Id} \times \fr} \ar[d]^{\sigma}  S'\cR\cF \\
\nt &\ar[l]_\fr \cR\cF
}
\end{equation}

\subsection{Cyclotomic spectra}\label{sectspectra}
The notion for spectra corresponding to the construction reviewed above is that  of {\em cyclotomic spectrum}.
The general notion of G-spectrum is here used in the special case G$=\T$, where the compact Lie group $\T$ has the special property that it is canonically isomorphic to its quotient $\T/C$ by any non-trivial closed subgroup $C\subset \T$. To an integer $n\in \nt$ corresponds the subgroup $C_n$ of $n$-th roots of unity;  there is an action $\Psi$ of the monoid $\nt$ on $\T$-spectra, by restriction to the fixed points of $C_n$ and then viewing the corresponding 
$\T/C_n$-spectrum as a $\T$-spectrum (see \cite{HM} Section 2, \cite{HM1} page 9).

\begin{definition}\label{defncyclospec}
A cyclotomic spectrum $E$ is  a $\T$-spectrum together with morphisms of $\T$-spectra $\epsilon_k:\Psi^k(E)\to E$ such that the following diagrams commute 
\begin{equation*}\label{tspec}
\xymatrix{
\Psi^k(\Psi^r(E)) \ar[d]_{\Psi^k(\epsilon_r)}  \ar[r]^-\simeq &\ar[d]^{\epsilon_{kr}}  \Psi^{kr}(E) \\
\Psi^k(E) \ar[r]^-{\epsilon_k} & E
}
\end{equation*}
\end{definition} 

By implementing \eqref{goodw1bis} to define the variant $S'\cR\cF$ of the category $S\cR\cF$, we obtain

\begin{lemma}
A cyclotomic spectrum $E$ defines a (covariant) functor from $S'\cR\cF$ to the category of RO$(\T)$-graded topological spaces.
\end{lemma} 
\proof Let consider the inclusion of  fixed points, \ie for $a=sb$ the map $(\theta,f_{1,s}):a\to b$.   Let $i_s^a$ be the inclusion  $E^{C_a}\subset E^{C_{a/s}}$. The action of $(\theta,f_{1,s})$ is $i^a_s\circ E(\frac \theta a)$. Thus (except for the inclusions) the action of $(\theta,f_{1,s}):a\to b$ composed with the action of   $(\tau, f_{1,q}):qa\to a$ is 
\[
 E\Big(\frac \theta a\Big)\circ  E\Big(\frac{\tau}{aq}\Big)=E\Big(\frac \theta a+\frac{\tau}{aq}\Big)=E\Big(\frac{q\theta +\tau}{aq}\Big)
\]
which agrees with \eqref{goodw1bis}. The  action of the maps  $f_{r,1}$  commutes with the action of $\T$ and corresponds to the composition of the map $\beta:E^{C_r}\to \Phi^{C_r}(E)$ with the structure map $\Phi^{C_r}(E)\to E$ of the cyclotomic spectrum. These  are all maps of $\T$-spectra by construction. Thus the parameters $r,q$ do not interfere with the validity of  \eqref{goodw1bis}. \endproof

\subsection{The pericyclic category $\Pi$}\label{sectpericyc}

There is a combinatorial version of the topological category $S'\cR\cF$: this is the pericyclic category $\Pi$ that we shall introduce here below. \newline
In the diagram \eqref{1o} we replace the monoid $(U(1)\ltimes \nt)^{\rm o}$ with the category $\lbto$. Indeed,  on the geometric realization of an epicyclic space in the sense of \cite{bu1} (\ie as a covariant functor from $\lbto$ to spaces) one has, by applying Theorem A of \opcit a natural {\em right}  action of the monoid $U(1)\ltimes \nt$, \ie  a left action of $(U(1)\ltimes \nt)^{\rm o}$. This shows that $\lbto$ is the combinatorial version of the monoid $(U(1)\ltimes \nt)^{\rm o}$.  In the following, we describe the missing entries (\ie the ?'s) in the pullback diagram  
\begin{equation}\label{3o}
\xymatrix{
 \lbto \ar[d]_{\bf{Mod}} & \ar[l]_{?} \ar[d]^{?}  ? \\
\nt &\ar[l]^\fr \cR\cF.
}\end{equation}

\begin{definition}\label{defnperi} The pericyclic category $\Pi$ is the subcategory of the product category $ \lbto\times\cR\cF$ sharing the same objects and whose morphisms 
\begin{equation}\label{22o}
\Hom_\Pi((n,a),(m,b))\subset \Hom_\lbto(n,m)\times \Hom_{\cR\cF}(a,b)
\end{equation}
are given by those pairs  $(h,f)\in \Hom_\lbto(n,m)\times \Hom_{\cR\cF}(a,b)$ fulfilling the equality
\begin{equation}\label{21o}
 \fr(f)={\rm Mod}(h).
\end{equation}
\end{definition}

 Thus an object of $\Pi$ is  a pair $(n,a)$ with $n\geq 0$, $n\in \N$ and $a\in \nt$. 
Then we complete  \eqref{3o} to a pullback diagram as follows
\begin{equation*}\label{4o}
\xymatrix{
 \lbto \ar[d]_{\bf{Mod}} & \ar[l]_{\lambda} \ar[d]^{\pi}  \Pi \\
\nt &\ar[l]^\fr {\cR\cF}
}\end{equation*}
by implementing the two projections arising from \eqref{22o}, namely the functors $\lambda:\Pi\longrightarrow \lbto$ and $\pi:\Pi\longrightarrow {\cR\cF}$ defined by
\[
\lambda(n,a):=n, \quad \lambda(h,f):=h,\qquad \pi(n,a):=a,\quad \pi(h,f):=f.
\]
To obtain a natural section of $\lambda$ we introduce the following further morphisms in $\Pi$.  For $n\geq 0$, $s>0$ and $d>0$ integers, with $m=s(n+1)-1$, one lets
\begin{equation*}\label{lift0}
\chi(n,s,d):=	(\id_n^{s,{\rm o}}, f_{1,s})\in \Hom_\Pi((n,sd),(m,d)).
\end{equation*}
This is meaningful since $(\id_n^{s,{\rm op}}, f_{1,s})\in \Hom_\lbto(n,m)\times \Hom_{\cR\cF}(sd,d)
$ fulfills \eqref{21o}.
 Let $W$ be the set of morphisms of the form $\id \times f_{r,1}:(n,a)\to (n,b)$ (for $br=a$). We localize the category $\Pi$ on $W$: this is the category $\Pi[W^{-1}]$ obtained by inverting the morphisms in $W$. One obtains in this way a small category equivalent to $\lbto$
 
\begin{lemma}\label{percycloc}
The  functor $\theta:\lbto\longrightarrow \Pi[W^{-1}]$ 
\[
\theta(n):=(n,1), \quad  \theta(\ell)=\ell \times f_{1,1},~\forall \ell\in \Hom_\Lambda(n,m),  \quad \theta(\id_n^{s,{\rm o}})=\chi(n,s,1)
\circ f_{s,1}^{-1} 
\]
defines an equivalence of categories.
\end{lemma}
\proof By construction, the restriction of $\theta$ to the cyclic category $\Lambda$ is an isomorphism with a subcategory of $\Pi$. Let $p_n^s:=\chi(n,s,1)
\circ f_{s,1}^{-1}:(n,1)\to (n,s)\to (s(n+1)-1,1)
$. One easily checks that these morphisms of $\Pi[W^{-1}]$ fulfill the presentation of $\lbto$ as an extension of the cyclic category. Indeed, this follows from the construction of $\Pi$ as a  subcategory of $\lbto \times \cR\cF$, so that one can check the relations in $\lbto \times \cR\cF$ using the fact that they hold in $\lbto$. This proves that $\theta$ is a functor. Moreover, the composite $\lambda\circ \theta$ is the identity thus $\theta$ is faithful. Finally, after inverting the morphisms in $W$ one sees  that an object $(n,a)$ of $\Pi$ is isomorphic to the object $(n,1)$ which is in the range of $\theta$. Thus $\theta$ is an equivalence of categories.\endproof

 Lemma \ref{percycloc} shows that covariant functors $\Pi[W^{-1}]\longrightarrow \Se$ correspond to epicyclic sets in the sense of \cite{bu1}. 
For an epicyclic space in the sense of \cite{DGM} (Definition 6.2.3.1), \ie of Definition \ref{defnepicycspace} we have

\begin{proposition}\label{bothepicyclic} Let $(Y,\phi)$ be an epicyclic space as in 
Definition \ref{defnepicycspace}. The following equalities define a covariant functor:
\begin{equation*}\label{beta0}  \beta :\Pi \longrightarrow \Se,\qquad \beta(n,a)={\text{ sd}}_a(Y)^{C_a}(n)\end{equation*}
\begin{equation}\label{beta1}\beta((\ell , \id_a))={\text{ sd}}_a(Y)^{C_a}(\ell):\beta(n,a)\to \beta(m,a), \ \forall \ell\in \Hom_\Lambda(n,m)\end{equation}
\begin{equation}\label{beta2}
\beta(\chi(n,u,v)):\beta(n,uv)\subset \beta(u(n+1)-1,v)
\end{equation}
\begin{equation*}\label{beta3}
\beta((\id , f_{q,1}))=\phi_q:\beta(n,qa)={\text{ sd}}_{aq}(Y)^{C_{aq}}(n)\to \beta(n,a)={\text{sd}}_a(Y)^{C_a}(n).
\end{equation*}
The functor $\beta$  encodes faithfully the epicyclic space $(Y,\phi)$. 
\end{proposition}
\proof Let us first explain the meaning of the inclusion  in \eqref{beta2}. By construction, one has for all $a>0$ and $n\geq 0$
\[
{\text{sd}}_a(Y)(n)=Y({\text{ sd}}^a(n))=Y(a(n+1)-1),
\]
so that 
\[
{\text{sd}}_{uv}(Y)(n)=Y(uv(n+1)-1)={\text{sd}}_{v}(Y)(u(n+1)-1).
\]
The fixed points under the cyclic subgroups involve the action of the subgroup $C_{uv}$ for $\beta(n,uv)={\text{ sd}}_{uv}(Y)^{C_{uv}}(n)$, and its subgroup $C_v\subset C_{uv}$ for $\beta(u(n+1)-1,v)={\text{ sd}}_{v}(Y)^{C_{v}}(u(n+1)-1)$. Thus the inclusion in \eqref{beta2} is the inclusion of fixed points. The structure maps $\phi_q$ are cyclic maps ${\text{sd}}_{aq}(Y)^{C_{aq}}\to {\text{sd}}_a(Y)^{C_a}$ and this shows that $\beta((\id , f_{q,1}))$ commutes with the maps $\beta((\ell , \id_a))$ of \eqref{beta1} defining the cyclic action. Moreover $\beta((\id , f_{q,1}))$ commutes with the inclusion of fixed points \eqref{beta2}, \ie the following diagram is commutative for any $n,q,u,v$
\begin{equation*}\label{14o}
\xymatrixcolsep{3pc}\xymatrix{
\beta(n,quv)\ar[d]_{\beta((\text{Id}, f_{q,1}))}  \ar[rr]^-{\beta{(\chi(n,u,qv))}}&& \ar[d]^{\beta((\text{Id} , f_{q,1}))}  \beta (u(n+1)-1,qv) \\
\beta(n,uv) \ar[rr]_-{\beta{(\chi(n,u,v))}}&& \  \  \beta (u(n+1)-1,v).
}
\end{equation*}
It remains to check that the inclusion of fixed points \eqref{beta2} fulfills the relation of $\id_n^{s,{\rm o}}$ in the category $\lbto$ with respect to the cyclic action \eqref{beta1}. This follows at the combinatorial level from the above discussion of the inclusion of fixed points. \endproof 
Proposition \ref{bothepicyclic} applies, in particular, to the cyclic nerve of a small category  $\cC$ (\cf~ \cite{DGM} Example 6.2.3.3). This is the cyclic set

\begin{equation*}\label{cyclicnerve}
B^{\rm cy}_n\cC=\{c_n\stackrel{f_0}{\leftarrow} c_0\stackrel{f_1}{\leftarrow} c_1\leftarrow \cdots \stackrel{f_n}{\leftarrow} c_{n}\}
\end{equation*}
formed by $n+1$ composable morphisms of   $\cC$.  The face maps are defined by composition, and the degeneracies by insertion of identities, while the cyclic structure is given by cyclic permutations. In \cite{DGM} (Example 6.2.3.3), one gets an epicyclic space (in the sense of \cite{DGM}), using  the isomorphisms
\begin{equation}\label{phik}
\phi_k:\left({\rm sd}_k B^{\rm cy}\cC\right)^{C_k}\simeq B^{\rm cy}\cC.
\end{equation}
Proposition \ref{bothepicyclic} thus defines a covariant functor $B^{\rm epi}\cC: \Pi\longrightarrow \Se$. Moreover, since the $\phi_k$ are invertible, Lemma \ref{percycloc} applies and one gets an epicyclic  set in the sense of \cite{bu1}. This functor is the same as the epicyclic set associated to the cyclic nerve in \opcit Indeed,  the latter is the functor $\lbto\longrightarrow\Se$ which extends the canonical cyclic structure of $B^{\rm cy}\cC$ by the additional   map $p_n^k:B^{\rm cy}_n\to B^{\rm cy}_{k(n+1)-1}$ associated to the morphism $\id_n^{k,{\rm o}}\in \Hom_\lbto(n,k(n+1)-1)$. The map $p_n^k$ is the $k$-fold amplification, mapping the  element $(c_n\stackrel{f_0}{\leftarrow} c_0\stackrel{f_1}{\leftarrow} c_1\leftarrow \cdots \stackrel{f_n}{\leftarrow} c_{n})$ of $B^{\rm cy}_n$ into its repetition $k$-times
\begin{equation*}\label{cyclicnerve1}
c_n\stackrel{f_0}{\leftarrow} c_0\stackrel{f_1}{\leftarrow} c_1\leftarrow \cdots \stackrel{f_n}{\leftarrow} c_{n}\stackrel{f_0}{\leftarrow} c_0\stackrel{f_1}{\leftarrow} c_1\leftarrow \cdots \cdots \leftarrow c_0\stackrel{f_1}{\leftarrow} c_1\leftarrow \cdots \stackrel{f_n}{\leftarrow} c_{n}
\end{equation*}
The map $p_n^k$ is the inverse of $\phi_k$ in \eqref{phik} and this corresponds exactly to the equality $\theta(\id_n^{s,{\rm op}})=\chi(n,s,1)
\circ f_{s,1}^{-1} $ of Lemma \ref{percycloc}.

\subsection{The small category $\cR$ and points of the topos $\widehat \cR$}\label{sectR}

To some extent, one may view the kernel of either one of the two  functors in \eqref{resfrob} as the category $\cR$ whose objects are the same as for $\cR\cF$, and where there exists a unique morphism $f(a,b)\in \Hom_\cR(a,b)$ if $b\vert a$, while otherwise there is no morphism.  Since the opposite category $\cR^{\rm o}$ has a morphism $a\to b$ exactly when $a\vert b$, it would seem that after a process of completion providing the points of the topos $\hat \cR$ as colimits of objects of $\cR^{\rm o}$, these limits ought to be naturally interpreted in terms of supernatural numbers. It is natural to think of the objects of $\cR\cF$ by associating to the object $n\in\nt$ the cyclic subgroup $C_n\subset \Q/\Z$ given by the $n$-torsion:
\[
C_n:=\{r\in \Q/\Z\mid n r=0\}=\frac 1n\Z/\Z.
\]
\begin{lemma}	
 \label{ptsD}The category of points of the topos $\widehat \cR$ is equivalent to the category of subgroups of $\Q$ containing $\Z$, with morphisms  given by inclusion.
\end{lemma}
\proof Let $F:\cR\longrightarrow\Se$ be a flat functor. Let $X_a$ be  the image by $F$ of the object $a$ of $\cR$. One has a unique map $F(a,b):X_{a}\to X_b$ for any $a,b\in \nt$ with $b\vert a$.  The filtering conditions 
for the category $\int_{\cR}\, F$ read as
\begin{enumerate}
\item $X_a\neq \emptyset$ for some $a\in \nt$.
\item For any $x\in X_a$, $y\in X_b$ $\exists$ $c\in \nt$ and $z\in X_c$  such that $a\vert c$, $b\vert c$ and $F(c,a)z=x$, $F(c,b)z=y$.
\end{enumerate}
The third filtering condition  is automatically fulfilled since for objects $i,j$ of  $\int_{\cR}\, F$ there is at most one morphism  $i  \to j$, thus $\alpha, \beta:i\to j$ implies $\alpha=\beta$. The second condition implies, by taking $a=b$, that there is at most one element in $X_a$ for any $a\in \nt$. Let $J:=\{a\in \nt \mid X_a\neq \emptyset\}$. One then has
\begin{equation}\label{pointsD}
a\in J, \quad b\vert a\implies b\in J, \quad  a,b\in J\implies \exists c\in J, ~ a\vert c, ~ b\vert c.
\end{equation}
The first implication follows from the existence of the morphism $F(a,b):X_{a}\to X_b$ while the second one follows from  $(2)$. Thus  $1\in J$ and the flat functors $F:\cR\longrightarrow\Se$ correspond to subsets $J\subset \nt$ fulfilling  \eqref{pointsD}. A morphism of functors $F\to F'$ exists if and only if $J\subset J'$, and in that case  is unique since both the sets involved have one element. We  use the following correspondence between subsets $J\subset \nt$ fulfilling  \eqref{pointsD} and subgroups $H$ such that $\Z\subset H\subset \Q$ 
\begin{equation*}\label{hj}
H^J:=\bigcup_{n\in J}\frac 1n\Z\subset \Q, \qquad  J_H:=\{n\in \nt\mid \frac 1n\in H\}.
\end{equation*}
The conditions \eqref{pointsD} show that $H^J$ is a filtering union of subgroups and is hence a subgroup of $\Q$ which contains $\Z$ since $1\in J$. Conversely, for $\Z\subset H\subset \Q$, the subset $J_H\subset \nt$ fulfills  the conditions \eqref{pointsD}. Finally, it is easy to see that the maps $J\mapsto H^J$ and $H\mapsto J_H$ are inverse of each other. \endproof

In \cite{CCas} we proved that  the non-trivial subgroups of $\Q$ are parameterized by the quotient space $\A_\Q^f/\hat\Z^*$ (where $\A_\Q^f$ are the finite ad\`eles of $\Q$), while the subgroups of $\Q$ containing $\Z$ are parameterized by the subset $\hat \Z/\hat\Z^*\subset \A_\Q^f/\hat\Z^*$. This latter subset surjects onto the quotient $\Q_+^\times\backslash\A_\Q^f/\hat\Z^*$ which parametrizes the points of $\widehat{\N^{\times}}$.
Here, one has a natural geometric morphism of toposes associated with the functor
\[
\rho:\cR\longrightarrow \nt, \qquad \rho(f(a,b)):=a/b
\]
The next proposition provides a topos theoretic interpretation  of the quotient map  $\A_\Q^f/\hat\Z^* \to \Q_+^\times\backslash\A_\Q^f/\hat\Z^*$.
\begin{proposition} \label{ptsD1}
The geometric morphism of toposes  $\widehat \cR\longrightarrow \widehat{\N^{\times}}$ associated with the functor $\rho$ maps the category of points of  $\hat \cR$ to the category of points of $\widehat{\N^{\times}}$ as the quotient map
\begin{equation}\label{pointsmap1}
\hat \Z/\hat\Z^*\subset \A_\Q^f/\hat\Z^* \to \Q_+^\times\backslash\A_\Q^f/\hat\Z^*.
\end{equation}
\end{proposition}
\proof
For a point $p$ of  $\hat \cR$, let $J\subset \nt$ be the associated subset. The pullback $p^*$  associates to a contravariant functor $Z:\cR\longrightarrow\Se$ the set 
\begin{equation}\label{pointsmap}
\left(\coprod_{a\in J}Z_a\right)/\sim ~=\varinjlim_{a\in J} Z_a
\end{equation}
as the colimit of the filtering diagram of  sets $Z_a\to Z_b$ for $a\vert b$ ($Z$ is contravariant). The contravariant functor $y(\bullet):\nt\to \Se$ which associates to the unique object $\bullet$ of $\nt$ the set $\nt=\Hom_{\nt}(\bullet,\bullet)$ on which $\nt$ acts by multiplication, once composed with $\rho$ defines  the following contravariant  functor \begin{equation}\label{pointsmap2}
Z=y(\bullet)\circ \rho:\cR\longrightarrow\Se, \qquad  Z(a)=\nt, \qquad Z(f(a,b))=y(\bullet)(a/b).
\end{equation}
Thus the flat functor $F:\nt\longrightarrow\Se$ associated with the image of the point $p$ by the geometric morphism, is given by the action of $\nt$ on the filtering colimit \eqref{pointsmap} for the functor $Z$ of \eqref{pointsmap2}. The computation of this colimit gives 
\[
\varinjlim_{a\in J}\nt=\bigcup_{a\in J}\frac 1a\Z_+=(H_J)_+.
\]
This shows that at the level of the subgroups of $\Q$ the map of points associates with $H$ ($\Z\subset H\subset \Q$) the group $H$ viewed as an abstract ordered group. This corresponds precisely to the map \eqref{pointsmap1}. \endproof


\begin{thebibliography}{99} 

\bibitem{bu1}  D.~Burghelea, Z.~Fiedorowicz, W.~Gajda,  {\em Power maps and epicyclic spaces}, Journal of Pure and Applied Algebra 96 (1994), 1--14.

\bibitem{CW} O.~Caramello, N.~Wentzlaff, {\em Cyclic theories}. Appl. Categ. Structures 25 (2017), no. 1, 105--126.

\bibitem{Ober81} A.~Connes, {\em Spectral sequence and homology
of currents for operator algebras}, {\it Math. Forschungsinst. Oberwolfach Tagungsber., 41/81,  Funktionalanalysis und $C^*$-Algebren, 27-9/3-10}, 1981.

\bibitem{CoExt} A.~Connes, {\em  Cohomologie cyclique et foncteurs
${\rm Ext}\sp n$}.  C. R. Acad. Sci. Paris S\'er. I Math. 296
(1983), no. 23, 953--958.

\bibitem{CoIHES} A.~Connes, {\em Noncommutative differential
geometry}.  Inst. Hautes \'Etudes Sci. Publ. Math.  No. 62 (1985), 257--360.


\bibitem{Co-book}  A.~Connes, {\it Noncommutative geometry}, Academic Press (1994).

\bibitem{Cad} A.~Connes, {\em Trace formula in noncommutative geometry and the zeros of the Riemann zeta function},  Selecta Math. (N.S.)  5  (1999),  no. 1, 29--106.

\bibitem{CC6} A.~Connes, C.~Consani, {\em Cyclic homology, Serre's local factors and the $\lambda$-operations},  J. K-Theory {\bf 14} (2014), no. 1, 1--45.

 
\bibitem{CCproj} A.~Connes, C.~Consani, {\em Projective geometry in characteristic one and the epicyclic category}, Nagoya Mathematical Journal  217 (2015), 95--132.  

\bibitem{topos} A.~Connes, C.~Consani, {\em  The cyclic and epicyclic sites}. Rend. Semin. Mat. Univ. Padova 134 (2015), 197--237.

\bibitem{cyctop} A.~Connes, C.~Consani, {\em Cyclic Structures and the Topos of Simplicial Sets},   J. Pure Appl. Algebra {\bf 219} (2015), no. 4, 1211--1235.

\bibitem{CCas} A.~Connes, C.~Consani, {\em Geometry of the Arithmetic Site}. Adv. Math. 291 (2016), 274--329.

\bibitem{CCprel} A.~Connes, C.~Consani, {\em Absolute algebra and Segal's Gamma sets},  J. Number Theory 162 (2016), 518--551.

\bibitem{CCscal1} A.~Connes, C.~Consani, \emph{Geometry of the Scaling Site}.  Selecta Math. (N.S.) \textbf{23} (2017), no. 3, 1803--1850.

\bibitem{CCscal3} A.~Connes, C.~Consani, {\em Homological Algebra in Characteristic One}, 
Higher Structures Journal {\bf 3} (2019), no. 1, 155--247.

\bibitem{CCgromov} A.~Connes, C.~Consani, {\em $\overline{\Spec\Z}$ and the Gromov norm}, Theory and Applications of Categories, {\bf 35}, No. 6, (2020), 155--178. 

\bibitem{CCAtiyah} A.~Connes, C.~Consani, {\em Segal's gamma rings and universal arithmetic}. Q. J. Math. 72 (2021), no. 1--2, 1--29.

\bibitem{schemeF1}   A.~Connes, C.~Consani, {\em On Absolute Algebraic Geometry, the affine case}, Adv. Math. 390 (2021), Paper No. 107909, 44p. 

\bibitem{RRinv}   A.~Connes, C.~Consani, {\em Riemann Roch for $\spzb$},  Preprint (2022).
Available at \url{https://arxiv.org/abs/2205.01391}


\bibitem{DGM} Dundas, Bjorn Ian; Goodwillie, Thomas G.; McCarthy, Randy {\em The local structure of algebraic K-theory}. Algebra and Applications, 18. Springer-Verlag London, Ltd., London, 2013.
\bibitem{FT} B.~Feigin, B.~Tsygan, {\em Additive K-theory}. K-theory, arithmetic and geometry (Moscow, 1984--1986), 67--209, Lecture Notes in Math., 1289, Springer, Berlin, 1987.

\bibitem{G} J. S. Golan {\em Semirings and their Applications}. Kluwer Academic Publishers, Dordrecht/Boston/London 1999.

\bibitem{Ober87} T.~Goodwillie, {\em A Topologist View of Cyclic Homology}, {\it Math. Forschungsinst. Oberwolfach Tagungsber., 26/87,  Zyklische Kohomologie und ihre Anwendungen, 14-6/20-6}, 1987.

\bibitem{GW} T.~Goodwillie, {\em Notes on the cyclotomic trace} Lecture notes for a series of seminar talks at MSRI, 1990-1991.

\bibitem{HM} L.~Hesselholt, I.~Madsen, {\em On the K-theory of finite algebras over Witt vectors of perfect fields}. Topology 36 (1997), no. 1, 29--101. 


\bibitem{HM1} L.~Hesselholt, I.~Madsen, {\em On the K-theory of local fields}, Annals of Mathematics, 158 (2003), 1--113.


\bibitem{Kaledin} D.~Kaledin, {\em Hochschild-Witt complex}, Advances in Mathematics 351 (2019) 33--95.

\bibitem{Loday} J.L.~Loday, {\em Cyclic homology}. Grundlehren der
Mathematischen Wissenschaften, 301. Springer-Verlag, Berlin, 1998.

\bibitem{Miyazaki} K.~Miyazaki, {\em Distinguished elements in a space of distributions}. J. Sci. Hiroshima Univ. Ser. A 24 (1960), 527--533.
\end{thebibliography}
\end{document}